\documentclass[12pt,reqno]{amsart}

\usepackage{amssymb, amsmath}
\usepackage{tikz-cd}

\usepackage{hyperref}

\numberwithin{equation}{section}

\newtheorem{theorem}{Theorem}[section]

\newtheorem{lemma}[theorem]{Lemma}

\theoremstyle{definition}
\newtheorem{definition}[theorem]{Definition}
\newtheorem{remark}[theorem]{Remark}
\newtheorem{notation}[theorem]{Notation}

\begin{document}

\baselineskip=15pt

\title[Genuinely ramified maps and stable vector bundles]{Genuinely ramified maps and stable vector bundles}

\author[I. Biswas]{Indranil Biswas}

\address{School of Mathematics, Tata Institute of Fundamental
Research, Homi Bhabha Road, Mumbai 400005, India}

\email{indranil@math.tifr.res.in}

\author[S. Das]{Soumyadip Das}

\address{School of Mathematics, Tata Institute of Fundamental
Research, Homi Bhabha Road, Mumbai 400005, India}

\email{dass@math.tifr.res.in}

\author[A. J. Parameswaran]{A. J. Parameswaran}

\address{School of Mathematics, Tata Institute of Fundamental
Research, Homi Bhabha Road, Mumbai 400005, India}

\email{param@math.tifr.res.in}



\date{}

\begin{abstract}
Let $f\, \colon \, X\, \longrightarrow\, Y$ be a separable finite surjective map between
irreducible normal projective varieties defined over an algebraically closed field, such that
the corresponding homomorphism between \'etale fundamental groups
$f_*\, \colon \, \pi_1^{\rm et}(X)\,\longrightarrow\, \pi_1^{\rm et}(Y)$ is surjective. Fix a
polarization on $Y$ and equip $X$ with the pullback, by $f$, of this polarization on $Y$.
Given a stable vector bundle $E$ on $X$, we prove that there is a vector bundle $W$ on $Y$ with
$f^*W$ isomorphic to $E$ if and only if the direct image $f_*E$ contains a stable vector bundle $F$ such that
$$
\frac{{\rm degree}(F)}{{\rm rank}(F)}\,=\, \frac{1}{{\rm degree}(f)}\cdot \frac{{\rm degree}(E)}{{\rm rank}(E)}\,.
$$
We also prove that $f^*V$ is stable for every stable vector bundle $V$ on $Y$.
\end{abstract}

\maketitle

\section{Introduction}

Here we continue, from \cite{PS}, \cite{BP1}, \cite{BP2}, the investigations of the direct 
image $f_*{\mathcal O}_X$, where $f\, :\, X\, \longrightarrow\, Y$ is a separable finite
surjective map between irreducible normal projective varieties defined over an 
algebraically closed field $k$. First we briefly recall the main results of \cite{BP1}, 
\cite{BP2}.

The above map $f$ is called genuinely ramified if the homomorphism between
\'etale fundamental groups
$$f_*\, \colon \, \pi_1^{\rm et}(X)\,\longrightarrow\, \pi_1^{\rm et}(Y)$$
induced by $f$ is surjective.
Consider the unique maximal pseudostable subbundle ${\mathbb S}\, \subset\, f_*{\mathcal O}_X$;
the existence of this subbundle is proved in \cite{BP2}.
The map $f$ is genuinely ramified if and only if ${\rm rank}({\mathbb S})\,=\, 1$ \cite{BP2}. When
$\dim X\,=\,1\,=\, \dim Y$, this was proved earlier in \cite{BP1}. When
$\dim X\,=\,1\,=\, \dim Y$, if $E$ is a stable vector bundle on $Y$, and $f$ is genuinely ramified, then $f^*E$ is 
also stable \cite{BP1}.

Here we prove the following (see Theorem \ref{thm1}):

\begin{theorem}\label{thmi1}
Let $f\, :\, X\, \longrightarrow\, Y$ be a finite separable
surjective map between irreducible normal projective varieties
defined over an algebraically closed field $k$. Then the following four statements
are equivalent.
\begin{enumerate}
\item The map $f$ is genuinely ramified.

\item There is no nontrivial \'etale covering $g\, :\, Y'\,\longrightarrow\, Y$ satisfying the condition
that there is a map $g'\, :\, X\,\longrightarrow\, Y'$ for which $g\circ g'\,=\, f$.

\item The fiber product $X\times_Y X$ is connected.

\item $\dim H^0(X,\, (f^*f_*{\mathcal O}_X)/(f^*f_*{\mathcal O}_X)_{\rm torsion})\,=\, 1$.
\end{enumerate}
\end{theorem}

This theorem leads to the following generalization, to higher dimensions, of the above mentioned result
for curves (see Theorem \ref{thm2}):

\begin{theorem}\label{thmi2}
Let $f\, :\, X\, \longrightarrow\, Y$ be a finite separable surjective map
between irreducible normal
projective varieties such that $f$ is genuinely ramified. Let $E$ be a stable vector bundle on
$Y$. Then the pulled back vector bundle $f^*E$ is also stable.
\end{theorem}

It should be noted that given a map $f$ as above, if $f^*E$ is a stable vector bundle on $Y$ for a vector bundle $E$ on $X$, then $E$ is stable.

Given a stable vector bundle $E$ on $X$ we may ask for a criterion for it to descend to $Y$ (meaning
a criterion for $E$ to be the pullback of a vector bundle on $Y$). The following criterion is deduced
using Theorem \ref{thmi2} (see Theorem \ref{thm3}):

\begin{theorem}\label{thmi3}
Let $f\, \colon \, X\, \longrightarrow\, Y$ be a genuinely ramified map of irreducible normal
projective varieties. Let $E$ be a stable vector bundle on $X$. Then there is a vector bundle $W$ on $Y$ with
$f^*W$ isomorphic to $E$ if and only if $f_*E$ contains a stable vector bundle $F$ such that
$$
\frac{{\rm degree}(F)}{{\rm rank}(F)}\,=\, \frac{1}{{\rm degree}(f)}\cdot \frac{{\rm degree}(E)}{{\rm rank}(E)}\,.
$$
\end{theorem}

Let $M$ be a smooth projective variety over $k$. If the \'etale fundamental group
of $M$ is trivial, then the stratified fundamental group of $M$ is also trivial \cite{EM}.
Let $f\, :\, M\, \longrightarrow\, {\mathbb P}^d_k$ be a finite
separable map, where $d\,=\, \dim M$.
The above mentioned result of \cite{EM} has the following equivalent reformulation:
If $f$ induces an isomorphism of the \'etale fundamental groups, then it
induces an isomorphism of the stratified fundamental groups. In view of this
reformulation, it is natural to ask the following question:

If $f\, :\, X\, \longrightarrow\, Y$ is a finite separable
surjective map between irreducible normal projective varieties such that the
corresponding homomorphism of the \'etale fundamental groups is surjective, then does
$f$ induce a surjection of the stratified fundamental groups?

We hope that Theorem \ref{thmi2} would shed some light on this question.

\section{Genuinely ramified maps}

Let $k$ be an algebraically closed field; there is no assumption on the characteristic of $k$.

Following \cite{PS}, \cite{BP1}, \cite{BHS} we define:

\begin{definition}\label{def1}
A separable finite surjective map
$$
f\,\, \colon \,\, X\, \longrightarrow\, Y
$$
between irreducible normal projective $k$-varieties is said to be \textit{genuinely ramified}
if the homomorphism between \'etale fundamental groups
$$f_*\, \colon \, \pi_1^{\rm et}(X)\,\longrightarrow\, \pi_1^{\rm et}(Y)$$
induced by $f$ is surjective.
\end{definition}

Let
\begin{equation}\label{e1}
f\, \colon \, X\, \longrightarrow\, Y
\end{equation}
be a separable finite surjective map between irreducible normal projective varieties.
Since $X$ and $Y$
are normal, and $f$ is a finite surjective map, the direct image $f_*{\mathcal O}_X$
is a reflexive sheaf on $Y$ whose rank coincides with the degree of the map $f$.

To define the degree of a torsionfree sheaf on $Y$, we fix a very ample
line bundle $L_Y$ on $Y$. Since $f$ is a finite map, the line bundle $f^*L_Y$ on $X$ is ample.
Equip $X$ with the polarization $f^*L_Y$. In what follows, the degree of the sheaves on
$Y$ and $X$ are always with respect to the polarizations $L_Y$ and $f^*L_Y$ respectively. For
a torsionfree sheaf $F$, define
$$
\mu(F)\,:=\, \frac{{\rm degree}(F)}{{\rm rank}(F)}\,\in\, \mathbb Q\, .
$$
We recall
that a torsionfree sheaf $E$ is called {\it stable} (respectively, {\it semistable}) if for all
subsheaves $V\, \subset\, E$, with $0\,<\, {\rm rank }(V)\,<\, {\rm rank}(E)$, the inequality
$$
\mu(V)\, <\, \mu(E) \ \ {\rm (respectively, }\ \mu(V)\, \leq\, \mu(E){\rm )}
$$
holds. If $E_1\, \subset\, E$ is the first nonzero term of the Harder-Narasimhan filtration of $E$, then
$$
\mu_{\rm max}(E)\,:=\, \mu(E_1)
$$
\cite{HL}. In particular, if $E$ is semistable, then $\mu_{\rm max}(E)\,=\, \mu(E)$.

For any coherent sheaves $E$ and $F$ on $X$ and $Y$ respectively, the projection formula gives
$$
f_*{\rm Hom} (f^* F, \, E)\,=\, f_*((f^*F^*)\otimes E)\,=\, F^*\otimes f_*E\,=\, {\rm Hom}(F,\, f_*E)\, .
$$
Since $f$ is a finite map, this gives the following:
\begin{equation}\label{eq_adjoint}
H^0(X,\, {\rm Hom} (f^* F, \, E))\,=\, H^0(Y,\, f_*{\rm Hom} (f^* F, \, E))\,=\, H^0(Y, \, {\rm Hom}(F, \, f_* E))\, .
\end{equation}
Setting $E\,=\, {\mathcal O}_X$
and $F\,=\,{\mathcal O}_Y$ in \eqref{eq_adjoint} we conclude that the identification $f^*{\mathcal O}_Y\,
\stackrel{\sim}{\longrightarrow}\, {\mathcal O}_X$ produces an inclusion homomorphism
\begin{equation}\label{e0}
{\mathcal O}_Y\, \hookrightarrow\, f_*{\mathcal O}_X\, .
\end{equation}
{}From \eqref{e0} it follows immediately that $\mu_{\rm max}(f_*{\mathcal O}_X)\, \geq\, 0$. In fact,
\begin{equation}\label{e2}
\mu_{\rm max}(f_*{\mathcal O}_X)\, =\, 0
\end{equation}
\cite[Lemma 2.2]{BP1}; although Lemma 2.2 of \cite{BP1} is stated under the assumption that
$\dim Y\,=\,1$, its proof works for all dimensions.

We recall from \cite{BP2} the definition of a pseudo-stable sheaf.

\begin{definition}
A {\it pseudo-stable sheaf} on $Y$ is a semistable sheaf $F$ on $Y$ such that $F$ admits a filtration of subsheaves
\begin{equation}\label{fi1}
0 = F_0 \, \subsetneq \, F_1 \, \subsetneq \, \cdots \, \subsetneq \, F_{n-1} \, \subsetneq \, F_n \, = \, F
\end{equation}
satisfying the condition that $F_i/F_{i-1}$ is a stable reflexive sheaf with $\mu(F_i/F_{i-1})\,=\, \mu(F)$ for
every $1 \,\leq\, i \,\leq \,n$.

A {\it pseudo-stable bundle} on $Y$ is a pseudo-stable sheaf $F$ as above such that
\begin{itemize}
\item $F$ is locally free, and

\item for each $1 \,\leq\, i \,\leq \,n$ the quotient $F_i/F_{i-1}$ in \eqref{fi1} is locally free.
\end{itemize}
\end{definition}

The first one among the above two conditions is actually implied by the second condition.

For the map $f$ in \eqref{e1}
consider the maximal destabilizing subsheaf $HN_1(f_* \mathcal{O}_X)$ of the reflexive coherent sheaf
$f_* \mathcal{O}_X$ on $Y$; in other words, $HN_1(f_* \mathcal{O}_X)$ is the first nonzero term in the
Harder-Narasimhan filtration of $f_* \mathcal{O}_X$. Recall from \eqref{e2} that
$$\mu(HN_1(f_* \mathcal{O}_X)) = \mu_{\rm max}(f_* \mathcal{O}_X) = 0\,.$$
Since $\mathcal{O}_Y \,\subset\, HN_1(f_* \mathcal{O}_X)$ (see \eqref{e0}) is a stable locally free
subsheaf of degree $0$, we conclude that \cite[Theorem 4.3]{BP2} applies to $f_* \mathcal{O}_X$. Let
\begin{equation}\label{e3}
{\mathbb S}\, \subset\, f_*{\mathcal O}_X
\end{equation}
be the unique maximal pseudo-stable bundle of degree zero such that $(f_*{\mathcal O}_X)/{\mathbb S}$
is torsionfree \cite[Theorem 4.3]{BP2}. Note that we have
\begin{equation}\label{t1}
{\mathcal O}_Y\, \subset\,{\mathbb S}.
\end{equation}

\begin{notation}
For any $k$-variety $Z$, and a coherent sheaf $F$ on $Z$, we write
$$
F_{\rm tor} \subset F
$$
for the torsion subsheaf of $F$. Also, define
\begin{equation}\label{en}
H^0_{\mathbb F}(W,\, F)\, :=\, H^0(W,\, F/F_{\rm tor}).
\end{equation}
\end{notation}

In particular,
\begin{equation}\label{e4}
H^0_{\mathbb F}(X,\, f^*f_*{\mathcal O}_X)\, = \,
H^0(X,\, (f^*f_*{\mathcal O}_X)/(f^*f_*{\mathcal O}_X)_{\rm tor})\, .
\end{equation}
Consider the Cartesian diagram
\begin{equation}\label{e5}
\begin{tikzcd}[row sep=large]
X\times_Y X \arrow[r, "p_2"] \arrow[dr, phantom, "\square"] \ar[d, "p_1"] & X \arrow[d, "f"]\\
X \arrow[r, "\widetilde{f}:=f"] & Y
\end{tikzcd}
\end{equation}
where $p_1$ and $p_2$ are the natural projections to the first and second factors
respectively. Since $f$ is a finite morphism, by the proper
base change theorem of quasi-coherent sheaves we have
\begin{equation}\label{e6}
\widetilde{f}^*f_*{\mathcal O}_X \,=\, f^*f_*{\mathcal O}_X\,=\, (p_1)_* {\mathcal O}_{X\times_Y X}
\end{equation}
(see \cite[Part 3, Ch.~59, \S~59.55, Lemma 59.55.3]{St}).
Pulling back the inclusion map in \eqref{e0} to $X$ we obtain
$$
{\mathcal O}_X \, \subset\, f^*f_*{\mathcal O}_X\, \cong \,(p_1)_* {\mathcal O}_{ X\times_Y X}.
$$
This implies that
\begin{equation}\label{e4b}
\dim H^0_{\mathbb F}(X,\, f^*f_*{\mathcal O}_X)\, \geq\, 1
\end{equation}
(see \eqref{e4}).

We can also describe $H^0_{\mathbb F}(X,\, f^*f_*{\mathcal O}_X)$ in terms of the cohomology of a sheaf on
the reduced scheme $(X\times_Y X)_{\rm red}$. Let
\begin{equation}\label{cz}
{\mathcal Z}\, \,:=\, \, (X\times_Y X)_{\rm red}\,\, \subset\,\, X\times_Y X
\end{equation}
be the reduced subscheme. So ${\mathcal O}_{\mathcal Z}$ is a quotient of ${\mathcal O}_{X\times_Y X}$.
Let
\begin{equation}\label{ecj}
j \, :\, \mathcal{Z} \, \longrightarrow \, X\times_Y X
\end{equation}
be the closed immersion. Let
\begin{equation}\label{pi}
p'_1,\, p'_2\, \colon \, {\mathcal Z}\,\longrightarrow\,X
\end{equation}
be the restrictions to $\mathcal Z$ of the projections $p_1,\, p_2$ respectively in \eqref{e5}. For $i\,=\, 1,\,
2$, we have $p_i \circ j = p'_i$, where $j$ is the map in \eqref{ecj}. Then
\begin{eqnarray}
(p'_{1})_* {\mathcal O}_{\mathcal Z} & \, = \, & (p_1)_* j_* \mathcal{O}_{\mathcal{Z}}
\nonumber\\
 & \, = \, & (p_1)_* (\mathcal{O}_{X\times_Y X}/(\mathcal{O}_{X\times_Y X})_{\rm tor})
\nonumber\\
 & \, = \, & ((p_1)_*\mathcal{O}_{X\times_Y X})/((p_1)_*\mathcal{O}_{X\times_Y X})_{\rm tor}.
\label{tq}
\end{eqnarray}
To explain the last equality in \eqref{tq}, note the following:
\begin{itemize}
\item If $\mathcal V$ is a torsionfree sheaf on $X\times_Y X$, then
$(p_1)_*\mathcal V$ is a torsionfree sheaf on $X$.

\item If $\mathcal V$ is a torsion sheaf on $X\times_Y X$, then
$(p_1)_*\mathcal V$ is a torsion sheaf on $X$.
\end{itemize}
Hence \eqref{tq} holds.

In view of \eqref{tq}, invoking the isomorphism in \eqref{e6} it is concluded that
$$
(p'_{1})_* {\mathcal O}_{\mathcal Z}\, = \, (f^*f_*{\mathcal O}_X)/(f^*f_*{\mathcal O}_X)_{\rm tor}\, .
$$
This implies that
\begin{equation}\label{eq_torsion_reduced}
H^0(X,\, (p'_{1})_* {\mathcal O}_{\mathcal Z})\,=\, H^0_{\mathbb F}(X,\, f^*f_*{\mathcal O}_X)
\,=\, H^0(X,\, (p'_{2})_* {\mathcal O}_{\mathcal Z})
\end{equation}
(the notation in \eqref{en} is used).

\begin{theorem}\label{thm1}
As in \eqref{e1}, let $f\, :\, X\, \longrightarrow\, Y$ be a finite separable surjective
map between irreducible normal projective varieties. Then the following five statements
are equivalent.
\begin{enumerate}
\item The map $f$ is genuinely ramified.

\item There is no nontrivial \'etale covering $g\, :\, Y'\,\longrightarrow\, Y$ satisfying the condition
that there is a map $g'\, :\, X\,\longrightarrow\, Y'$ for which $g\circ g'\,=\, f$.

\item The fiber product $X\times_Y X$ is connected.

\item $\dim H^0_{\mathbb F}(X,\, f^*f_*{\mathcal O}_X)\,=\, 1$ (see \eqref{e4}).

\item The subsheaf ${\mathbb S}\, \subset\, f_*{\mathcal O}_X$ in \eqref{e3} coincides with the subsheaf
${\mathcal O}_Y$ in \eqref{e0};
in other words, the inclusion in \eqref{t1} is an equality. This is equivalent to the condition that
${\rm rank}({\mathbb S})\,=\, 1$.
\end{enumerate}
\end{theorem}

\begin{proof}
We will show that $(1) \Leftrightarrow (2)$, $(1) \Leftrightarrow (5)$, $(3) \Leftrightarrow (4)$, $(3) \Rightarrow 
(2)$ and $(5) \Rightarrow (4)$.

Let
$$f_*\, :\, \pi_1^{\rm et}(X)\,\longrightarrow\, \pi_1^{\rm et}(Y)$$
be the homomorphism between the \'etale fundamental groups induced by the map $f$.
Since the map $f$ is dominant, $f_*(\pi_1^{\rm et}(X))$ is a subgroup of
$\pi_1^{\rm et}(Y)$ of finite index. In fact, the index of this
subgroup is at most $\text{degree}(f)$. Let
$$
g\, :\, Y'\,\longrightarrow\, Y
$$ be the \'etale covering for this subgroup $f_*(\pi_1^{\rm et}(X))\, \subset\,
\pi_1^{\rm et}(Y)$. Then there is a morphism
$g'\, :\, X\,\longrightarrow\, Y'$ such that $g\circ g'\,=\, f$. To explain this map $g'$, fix a point
$y_0\, \in\, Y$ and also fix points $x_0\, \in\, X$ and $y_1\, \in\, Y'$ over $y_0$. Let $X'$
be the connected component of the fiber product $X\times_Y Y'$ containing the point $(x_0,\, y_1)$. The
given condition that $f_*(\pi_1^{\rm et}(X,\, x_0))\,=\, g_*(\pi_1^{\rm et}(Y',\, y_1))$ implies that the
\'etale covering $f'\, :\, X'\, \longrightarrow\, X$, where $f'$ is
the natural projection, induces an isomorphism of \'etale fundamental groups.
Hence $f'$ is actually an isomorphism. If $g'$ is the composition of maps $$X\, \stackrel{\sim}{\longrightarrow}\, X'
\, \longrightarrow\, Y'\, ,$$ where $X'\, \longrightarrow\, Y'$ is the natural projection, then
clearly, $g\circ g'\,=\, f$.

Therefore, the statement (2) implies the statement (1). Conversely, if there is a nontrivial \'etale covering
$g\, :\, Y'\,\longrightarrow\, Y$ and a morphism $g'\, :\, X\,\longrightarrow\, Y'$ such that
$g\circ g'\,=\, f$, then the homomorphism $f_*$ is not surjective, because the homomorphism $g_*$ is
not surjective. So the statements (1) and (2) are equivalent.

It is known that the statements (1) and (5) are equivalent; see \cite[Proposition 1.3]{BP2}.

We will now prove that the statements (3) and (4) are equivalent.

The fiber product $X\times_Y X$ is connected if and only if $\dim H^0({\mathcal Z},\, {\mathcal O}_{\mathcal Z})
\,=\,1$, where $\mathcal Z$ is constructed in \eqref{cz}. Therefore, $X\times_Y X$ is connected if and only if
$$\dim H^0(X,\, (p'_{1})_* {\mathcal O}_{\mathcal Z})\,=\, 1,$$ where $p'_1$ is the map in 
\eqref{pi}. Consequently, from \eqref{eq_torsion_reduced} it follows
immediately that the statements (3) and (4) are equivalent.

Next we will show that the statement (3) implies the statement (2).

Assume that (2) does not hold. So there is a nontrivial \'etale covering
$$g\, :\, Y'\,\longrightarrow\, Y$$ and a map
$g'\, :\, X\,\longrightarrow\, Y'$ such that $g\circ g'\,=\, f$. Since $g$ is \'etale, the fiber product
$Y'\times_Y Y'$ is a reduced normal projective scheme. However, $Y'\times_Y Y'$ is not connected. In fact,
the image of the diagonal embedding $Y'\, \hookrightarrow\, Y'\times_Y Y'$ is a connected component of
$Y'\times_Y Y'$. Note that this diagonal embedding is not surjective because $g$ is
nontrivial. Since the map
$$
(g',\, g')\, :\, X\times_Y X\,\longrightarrow\, Y'\times_Y Y'
$$
is surjective, and $Y'\times_Y Y'$ is not connected, we conclude that $X\times_Y X$ is
not connected. Hence the statement (3) implies the statement (2).

Finally, we will prove that the statement (5) implies the statement (4).

Assume that the statement (4) does not hold. So from \eqref{e4b} it follows that
\begin{equation}\label{em}
n\, :=\, \dim H^0_{\mathbb F}(X,\, f^*f_*{\mathcal O}_X)\,\geq \, 2\, .
\end{equation}

Choose a finite Galois field extension $k(Y)\,\subset\, k(X)\,\subset\,L$. Let $$\Gamma \, = \, \text{Gal}(L/k(Y))$$
be the corresponding Galois group. Let $M$ be the normalization of $X$ in $L$. So $M$ is an irreducible normal projective
variety with $k(M)\,=\,L$, and
$$
\widehat{f} \, \,\colon\,\, M \, \longrightarrow \, Y
$$
is a Galois (possibly branched) cover dominating the map $f$;
the Galois group for $\widehat{f}$ is $\text{Gal}(L/k(Y))\,=\, \Gamma$. Let
\begin{equation}\label{eg}
g \,\, \colon \,\, M \, \longrightarrow \, X
\end{equation}
be the morphism induced by the inclusion map $k(X)\,\hookrightarrow\,L$, so $f\circ g\,=\, \widehat{f}$.

Let
\begin{equation}\label{g1}
\varphi\, :\, {\mathcal O}_X\otimes H^0_{\mathbb F}(X,\, f^*f_*{\mathcal O}_X)\, \longrightarrow\,
(f^*f_*{\mathcal O}_X)/(f^*f_*{\mathcal O}_X)_{\rm tor}
\end{equation}
be the evaluation homomorphism. We will show that this homomorphism $\varphi$ of coherent sheaves
is injective. To prove this, first note that for a semistable sheaf $V$ on $Y$,
the pullback $(f^*V)/(f^*V)_{\rm tor}$ is semistable (see \cite[Remark 2.1]{BP1}; the proof in
\cite[Remark 2.1]{BP1} works for all dimensions). Therefore, from \eqref{e2} we conclude that
\begin{equation}\label{g2}
\mu_{\rm max}((f^*f_*{\mathcal O}_X)/(f^*f_*{\mathcal O}_X)_{\rm tor})\, =\, 0\, .
\end{equation}
Any coherent subsheaf ${\mathcal W}\, \subset\, {\mathcal O}_X\otimes H^0_{\mathbb F}
(X,\, f^*f_*{\mathcal O}_X)$
such that
\begin{enumerate}
\item $\text{degree}({\mathcal W})\,=\, 0$, and

\item the quotient $({\mathcal O}_X\otimes H^0_{\mathbb F}(X,\, f^*f_*{\mathcal O}_X))/{\mathcal W}$ is torsionfree
\end{enumerate}
must be of the form $${\mathcal O}_X\otimes W\, \subset\,{\mathcal O}_X\otimes H^0_{\mathbb F}(X,\,
f^*f_*{\mathcal O}_X)$$
for some subspace $$W\, \subset\, H^0_{\mathbb F}(X,\, f^*f_*{\mathcal O}_X).$$ Indeed, this follows from the
fact that $\text{degree}({\mathcal W})$ is a nonzero multiple of the degree of the pullback of the tautological bundle
by the rational map from $X$ to the Grassmannian corresponding to ${\mathcal W}$; if ${\rm rank}({\mathcal W})\,=\, a$, then
${\mathcal W}$ is given by a rational map from $X$ to the Grassmannian parametrizing $a$ dimensional subspaces of $k^{\oplus n}$ (this map is defined on the open subset over which $\mathcal W$ is a subbundle of ${\mathcal O}_X\otimes H^0_{\mathbb F}(X,\, f^*f_*{\mathcal O}_X)$). Therefore, from the given two conditions on ${\mathcal W}$ it follows that it is of the form ${\mathcal O}_X\otimes W$ for some subspace
$W\, \subset\, H^0_{\mathbb F}(X,\, f^*f_*{\mathcal O}_X)$. For the homomorphism $\varphi$ in
\eqref{g1} if ${\rm kernel}(\varphi)\, \not=\, 0$, then $\text{degree}({\rm kernel}(\varphi))\,=\, 0$; this follows from \eqref{g2} and the fact that ${\rm degree}({\mathcal O}_X\otimes H^0_{\mathbb F}(X,\, f^*f_*{\mathcal O}_X)) \,=\,0$ (these two together imply that the image of $\varphi$ lies in the
maximal semistable subsheaf of $(f^*f_*{\mathcal O}_X)/(f^*f_*{\mathcal O}_X)_{\rm tor}$).
Also, the quotient $({\mathcal O}_X\otimes H^0_{\mathbb F}(X,\, f^*f_*{\mathcal O}_X))/{\rm kernel}(\varphi)$ is torsionfree, because it is contained in $(f^*f_*{\mathcal O}_X)/(f^*f_*{\mathcal O}_X)_{\rm tor}$. Therefore, from the earlier observation we conclude that ${\rm kernel}(\varphi)$
is of the form ${\mathcal O}_X\otimes W\, \subset\,{\mathcal O}_X\otimes H^0_{\mathbb F}(X,\, f^*f_*{\mathcal O}_X)$ for some subspace $W\, \subset\, H^0_{\mathbb F}(X,\, f^*f_*{\mathcal O}_X)$. From this it follows that the evaluation map on any $w\, \in\, W\, \subset\, H^0_{\mathbb F}(X,\, f^*f_*{\mathcal O}_X)$ is identically zero. But this implies that $w\,=\, 0$. Therefore, we conclude that
$W\,=\, 0$, and the homomorphism $\varphi$ in \eqref{g1} is injective.

Let
\begin{equation}\label{es}
{\mathcal S}\,:=\, {\rm image}(\varphi)\, \subset\, (f^*f_*{\mathcal O}_X)/(f^*f_*{\mathcal O}_X)_{\rm tor}
\end{equation}
be the image of $\varphi$. Since $\varphi$ is injective, we conclude that
${\mathcal S}$ is isomorphic to ${\mathcal O}^{\oplus n}_X$ (see \eqref{em}). The pullback $g^*{\mathcal S}$
by the map $g$ in \eqref{eg} is free (equivalently, it is a trivial bundle)
because ${\mathcal S}$ is so. The pulled back homomorphism
$$
g^*\varphi\, :\, g^*{\mathcal S}\, \longrightarrow\,
(g^*f^*f_*{\mathcal O}_X)/(g^*f^*f_*{\mathcal O}_X)_{\rm tor}\,=\,
(\widehat{f}^*f_*{\mathcal O}_X)/(\widehat{f}^*f_*{\mathcal O}_X)_{\rm tor}
$$
is injective because $\varphi$ is injective and $\mathcal S$ is free (or
the bundle is trivial). Using $g^*\varphi$ we would consider $g^*{\mathcal S}$
as a subsheaf of $(\widehat{f}^*f_*{\mathcal O}_X)/(\widehat{f}^*f_*{\mathcal O}_X)_{\rm tor}$.

The Galois group $\Gamma\,=\, {\rm Gal}(\widehat{f})$ has a natural action on $\widehat{f}^*f_*{\mathcal O}_X$
because this sheaf is pulled back from $Y$. This action of $\Gamma$ on
$\widehat{f}^*f_*{\mathcal O}_X$ produces an action of $\Gamma$ on the torsionfree quotient
$(\widehat{f}^*f_*{\mathcal O}_X)/(\widehat{f}^*f_*{\mathcal O}_X)_{\rm tor}$.

Consider ${\mathcal Z}$ defined in \eqref{cz} together with the projections $p'_1$ and $p'_2$ from it in
\eqref{pi}. Both $(p'_1)^*{\mathcal S}$ and $(p'_2)^*{\mathcal S}$ are free because $\mathcal S$
is so. The natural isomorphism
$$
(p'_1)^*f^*f_*{\mathcal O}_X \, \stackrel{\sim}{\longrightarrow}\, (p'_2)^*f^*f_*{\mathcal O}_X,
$$
given by the fact that they are the pull backs of a single sheaf on $Y$, produces an isomorphism
\begin{equation}\label{it2}
((p'_1)^*f^*f_*{\mathcal O}_X)/((p'_1)^*f^*f_*{\mathcal O}_X)_{\rm tor} \, \stackrel{\sim}{\longrightarrow}\,
((p'_2)^*f^*f_*{\mathcal O}_X)/((p'_2)^*f^*f_*{\mathcal O}_X)_{\rm tor}\,.
\end{equation}
It can be shown that the isomorphism in \eqref{it2} takes the subsheaf
$$
(p'_1)^*{\mathcal S}\,\subset\, ((p'_1)^*f^*f_*
{\mathcal O}_X)/((p'_1)^*f^*f_*{\mathcal O}_X)_{\rm tor},
$$
where $\mathcal S$ is defined in \eqref{es}, to the subsheaf
$$
(p'_2)^*{\mathcal S}\, \subset\,((p'_2)^*f^*f_*{\mathcal O}_X)/((p'_2)^*f^*f_*{\mathcal O}_X)_{\rm tor}.
$$
Indeed, we have
$$
H^0(X,\, (p'_{1})_* {\mathcal O}_{\mathcal Z})\,=\, H^0({\mathcal Z},\, {\mathcal O}_{\mathcal Z})\,=\,
H^0(X,\, (p'_{2})_* {\mathcal O}_{\mathcal Z})
$$
(because $p'_1$ and $p'_2$ are finite maps), and therefore from \eqref{eq_torsion_reduced} it follows that both
$(p'_1)^*{\mathcal S}$ and $(p'_2)^*{\mathcal S}$ coincide with 
${\mathcal O}_{\mathcal Z}\otimes H^0({\mathcal Z},\, {\mathcal O}_{\mathcal Z})$.
The isomorphism in \eqref{it2} takes the subsheaf $(p'_1)^*{\mathcal S}$ to $(p'_2)^*{\mathcal S}$,
and the homomorphism $(p'_1)^*{\mathcal S}\,{\longrightarrow}\,(p'_2)^*{\mathcal S}$ given by the
isomorphism in \eqref{it2} coincides with the following automorphism of ${\mathcal O}_{\mathcal Z}
\otimes H^0({\mathcal Z},\, {\mathcal O}_{\mathcal Z})$ (after we
identify $(p'_1)^*{\mathcal S}$ and $(p'_2)^*{\mathcal S}$ with
${\mathcal O}_{\mathcal Z}\otimes H^0({\mathcal Z},\, {\mathcal O}_{\mathcal Z})$): The involution
$(x,\, y)\, \longmapsto\, (y,\, x)$
of $X\times_Y X$ produces an involution of ${\mathcal Z}\, =\, (X\times_Y X)_{\rm red}$. This
involution of ${\mathcal Z}$ in turn produces an involution 
$$
\delta\, :\, H^0({\mathcal Z},\, {\mathcal O}_{\mathcal Z}) \,{\longrightarrow}\,
H^0({\mathcal Z},\, {\mathcal O}_{\mathcal Z})\, .
$$
Define the automorphism
$$
{\rm Id}\otimes\delta\, :\, {\mathcal O}_{\mathcal Z}\otimes H^0({\mathcal Z},\, {\mathcal O}_{\mathcal Z})
\, \longrightarrow\, {\mathcal O}_{\mathcal Z}\otimes H^0({\mathcal Z},\, {\mathcal O}_{\mathcal Z})\, .
$$
The homomorphism $(p'_1)^*{\mathcal S}\,{\longrightarrow}\,(p'_2)^*{\mathcal S}$ given by the
isomorphism in \eqref{it2} coincides with the automorphism ${\rm Id}\otimes\delta$ of ${\mathcal O}_{\mathcal Z}
\otimes H^0({\mathcal Z},\, {\mathcal O}_{\mathcal Z})$, after we
identify $(p'_1)^*{\mathcal S}$ and $(p'_2)^*{\mathcal S}$ with
${\mathcal O}_{\mathcal Z}\otimes H^0({\mathcal Z},\, {\mathcal O}_{\mathcal Z})$.

Using the above observation that the isomorphism in \eqref{it2} takes the subsheaf $(p'_1)^*{\mathcal S}$
to $(p'_2)^*{\mathcal S}$ we will now show that the natural action
of the Galois group $\Gamma\,=\, {\rm Gal}(\widehat{f})$ on $(\widehat{f}^*f_*{\mathcal O}_X)/
(\widehat{f}^*f_*{\mathcal O}_X)_{\rm tor}$ preserves the subsheaf $g^*{\mathcal S}$, where $g$
is the map in \eqref{eg}; it was noted earlier
that $\Gamma$ acts on $(\widehat{f}^*f_*{\mathcal O}_X)/(\widehat{f}^*f_*{\mathcal O}_X)_{\rm tor}$.

To prove that the action of $\Gamma$ on $(\widehat{f}^*f_*{\mathcal O}_X)/
(\widehat{f}^*f_*{\mathcal O}_X)_{\rm tor}$ preserves $g^*{\mathcal S}$, let $Y_0\, \subset\,
Y$ be the largest open subset such that
\begin{itemize}
\item the map $f$ is flat over $Y_0$, and

\item the map $g$ is flat over $f^{-1}(Y_0)$.
\end{itemize}
Define $X_0\, :=\, f^{-1}(Y_0)$ and $M_0\,:=\, \widehat{f}^{-1}(X_0)$. The conditions on $Y_0$
imply that the restriction of $\widehat f$ to
$M_0$ is flat \cite[Ch~III, \S~9]{Ha}. In view of the descent criterion for sheaves under flat
morphisms, \cite{Gr}, \cite{Vi},
the above observation, that the isomorphism in \eqref{it2} takes the subsheaf
$(p'_1)^*{\mathcal S}$ to the subsheaf $(p'_2)^*{\mathcal S}$, implies that the restriction
${\mathcal S}\big\vert_{X_0}$ descends to a subsheaf of $(f_*{\mathcal O}_X)\big\vert_{Y_0}$.
Consequently, the action of $\Gamma$ on $(\widehat{f}^*f_*{\mathcal O}_X)\big\vert_{M_0}$ preserves
$(g^*{\mathcal S})\big\vert_{M_0}$ (as it is the pullback of a
sheaf on $Y_0$); note that $\widehat{f}^*f_*{\mathcal O}_X$ is torsionfree over $M_0$ because
the map $\widehat f$ is flat over $M_0$. The codimension of the complement $M\setminus M_0$ is at least
two. Therefore, given two vector bundle $A$ and $B$ on $M$ together with an isomorphism
$${\mathcal I}\, :\, A\big\vert_{M_0} \, \stackrel{\sim}{\longrightarrow}\, B\big\vert_{M_0}$$ over $M_0$, there
is a unique isomorphism
$$\widetilde{\mathcal I}\, :\, A \, \stackrel{\sim}{\longrightarrow}\, B$$
that restricts to ${\mathcal I}$ on $M_0$; recall that $M$ is normal. For any $\gamma\, \in\, \Gamma$,
set $g^*{\mathcal S}\,=\, A\,=\, B$ and ${\mathcal I}$ to be the action of $\gamma$ on
$(g^*{\mathcal S})\big\vert_{M_0}$; recall that $g^*{\mathcal S}$ is locally free and
the action of $\Gamma$ on $(\widehat{f}^*f_*{\mathcal O}_X)\big\vert_{M_0}$ preserves
$(g^*{\mathcal S})\big\vert_{M_0}$. Now from the above observation that $\mathcal I$ extends uniquely
to $M$ we conclude that the action of $\gamma$ on $(\widehat{f}^*f_*{\mathcal O}_X)
/(\widehat{f}^*f_*{\mathcal O}_X)_{\rm tor}$ preserves the subsheaf $g^*{\mathcal S}$.

{}From the above observation that $g^*{\mathcal S}$ is preserved by the natural action of $\Gamma$ on 
$(\widehat{f}^*f_*{\mathcal O}_X)/(\widehat{f}^*f_*{\mathcal O}_X)_{\rm tor}$ we will now
deduce that there is a locally free subsheaf
$$
V\, \subset\, f_*{\mathcal O}_X
$$
such that the image of $\widehat{f}^*V$ in $(\widehat{f}^*f_*{\mathcal O}_X)/(\widehat{f}^*f_*{\mathcal O}_X)_{\rm 
tor}$ coincides with the image of $g^*{\mathcal S}$ in $(\widehat{f}^*f_*{\mathcal O}_X)/(\widehat{f}^*f_*{\mathcal 
O}_X)_{\rm tor}$.

To prove this, consider the free sheaf
$$
{\mathcal O}_M\otimes H^0(M,\, g^*{\mathcal S})\, \longrightarrow\, M\, ;
$$
it corresponds to the trivial vector bundle on $M$ with fiber $H^0(M,\, g^*{\mathcal S})$. Let
\begin{equation}\label{ei}
\Phi\, :\, {\mathcal O}_M\otimes H^0(M,\, g^*{\mathcal S})\, \longrightarrow\, g^*{\mathcal S}
\end{equation}
be the evaluation map. This $\Phi$ is an isomorphism, because $g^*{\mathcal S}$ is free, or equivalently,
it is a trivial vector bundle (recall that ${\mathcal S}$ is free).
The action of $\Gamma$ on $g^*{\mathcal S}$ produces an action of $\Gamma$ on $H^0(M,\,
g^*{\mathcal S})$ which, coupled with the action of $\Gamma$ on $M$, produces an action of $\Gamma$ on
${\mathcal O}_M\otimes H^0(M,\, g^*{\mathcal S})$. The isomorphism $\Phi$
in \eqref{ei} is evidently $\Gamma$--equivariant for the actions of $\Gamma$.

Let
\begin{equation}\label{ga0}
\Gamma_0\, \subset\, \Gamma
\end{equation}
be the normal subgroup that acts trivially on the vector space $H^0(M,\, g^*{\mathcal S})$. Now using
$\Phi$ we conclude that $g^*{\mathcal S}$ descends to a vector bundle on the quotient $M/\Gamma_0$. Let
\begin{equation}\label{ga1}
{\mathcal S}_1\, \longrightarrow\, M/\Gamma_0
\end{equation}
be the descent of $g^*{\mathcal S}$. In other words, the pullback of ${\mathcal S}_1$ to $M$ has a
$\Gamma_0$--equivariant isomorphism with $g^*{\mathcal S}$. Note that ${\mathcal O}_M\otimes H^0(M,\,
g^*{\mathcal S})$ descends to ${\mathcal O}_{M/\Gamma_0}\otimes H^0(M,\, g^*{\mathcal S})$ on $M/\Gamma_0$,
because $\Gamma_0$ acts trivially on $H^0(M,\, g^*{\mathcal S})$, and
hence ${\mathcal S}_1$ is a trivial vector bundle.

For the action of $\Gamma$ on $M$, the isotropy subgroup $\Gamma_m\, \subset\, \Gamma$ any point $m\,\in\, M$
is actually contained in the subgroup $\Gamma_0$ in \eqref{ga0}. Indeed, this follows immediately from the fact that
$\Gamma_m$ acts trivially on the fiber of $(\widehat{f}^*f_*{\mathcal O}_X)/(\widehat{f}^*f_*{\mathcal O}_X)_{\rm tor}$
over $m$; recall that $\mathcal S$ is the image of the injective map $\varphi$ from a trivial vector bundle
(see \eqref{g1} and \eqref{es}). Consequently, the natural map
$$
M/\Gamma_0\, \longrightarrow\, Y
$$
given by $\widehat{f}$ is actually \'etale Galois with Galois group $\Gamma/\Gamma_0$. The
action of $\Gamma$ on $g^*{\mathcal S}$ produces an action of $\Gamma/\Gamma_0$ on ${\mathcal S}_1$ in \eqref{ga1};
it is a lift of the action of $\Gamma/\Gamma_0$ on $M/\Gamma_0$. Hence ${\mathcal S}_1$ descends to
a vector bundle ${\mathcal S}_2$ on $Y$.

We have ${\mathcal S}_2\, \subset\, {\mathbb S}$, where ${\mathbb S}$ is constructed in \eqref{e3}, and
also $${\rm rank}({\mathcal S}_2)\,=\, {\rm rank}({\mathcal S})\,=\, n\, \geq\, 2$$ (see \eqref{em}).
Hence the statement (5) in the theorem fails. Therefore, 
the statement (5) implies the statement (4). This completes the proof.
\end{proof}

\begin{theorem}\label{thm2}
Let $f\, :\, X\, \longrightarrow\, Y$ be a finite separable surjective map between irreducible normal
projective varieties such that $f$ is genuinely ramified. Let $E$ be a stable vector bundle on
$Y$. Then the pulled back vector bundle $f^*E$ is also stable.
\end{theorem}

\begin{proof}
In view of Theorem \ref{thm1}, the proof is exactly identical to the proof of Theorem 1.1 
of \cite{BP1}. The only point to note that in Proposition 3.5 of \cite{BP1}, the sheaves 
${\mathcal L}_i$ are no longer locally free. Now they are sheaves properly contained in 
${\mathcal O}_X$. But this does not affect the arguments. We omit the details.
\end{proof}

\section{A characterization of stable pullbacks}

\begin{lemma}\label{lem3}
Let $f\, \colon \, X\, \longrightarrow\, Y$ be a separable finite surjective map
of irreducible normal
projective varieties. Let $E$ be a semistable vector bundle on $X$. Then
$$
\mu_{\rm max}(f_*E)\, \leq\, \frac{1}{{\rm degree}(f)}\cdot\mu(E)\, .
$$
\end{lemma}

\begin{proof}
Let $F\, \subset\, f_*E$ be the first nonzero term of the Harder-Narasimhan filtration of $f_*E$,
so $$\mu_1\, :=\, \mu_{\rm max}(f_*E)\,=\, \mu(F).$$ Therefore,
$(f^*F)/(f^*F)_{\rm tor}$ is semistable \cite[Remark 2.1]{BP1}, and
$$
\mu((f^*F)/(f^*F)_{\rm tor})\,=\, {\rm degree}(f)\cdot\mu_1\,.
$$
The inclusion map $F\, \hookrightarrow\, f_*E$ gives a homomorphism $f^*F\, \longrightarrow\, E$
(see \eqref{eq_adjoint}), which in turn produces a homomorphism
$$
\beta\, :\, (f^*F)/(f^*F)_{\rm tor}\, \longrightarrow\, E.
$$
Since $\beta\, \not=\, 0$, and both $E$ and $(f^*F)/(f^*F)_{\rm tor}$ are
semistable, we have
$$
\mu(E)\, \geq\, \mu({\rm image}(\beta))
\, \geq\, \mu((f^*F)/(f^*F)_{\rm tor})\,=\, {\rm degree}(f)\cdot\mu_1.
$$
This proves the lemma.
\end{proof}

\begin{theorem}\label{thm3}
Let $f\, \colon \, X\, \longrightarrow\, Y$ be a genuinely ramified map of irreducible normal
projective varieties. Let $E$ be a stable vector bundle on $X$. Then there is a vector
bundle $W$ on $Y$ with
$f^*W$ isomorphic to $E$ if and only if $f_*E$ contains a stable vector bundle $F$ such that
$$
\mu(F)\,=\, \frac{1}{{\rm degree}(f)}\cdot \mu(E)\,.
$$
\end{theorem}

\begin{proof}
First assume that there is a vector bundle $W$ on $Y$ such that $f^*W$ is isomorphic to $E$.
It can be shown that $W$ is stable. Indeed, if
$$
S\, \subset\, W
$$
is a nonzero coherent subsheaf such that ${\rm rank}(W/S)\, >\, 0$
and $\mu(S)\, \geq\, \mu(W)$, then we have ${\rm rank}((f^*W)/(f^*S))\, >\, 0$
and $$\mu(f^*S)\,=\, \mu(S)\cdot {\rm degree}(f)\, \geq\, \mu(W)\cdot {\rm degree}(f)
\,=\, \mu(f^*W),$$
which contradicts the stability condition for $f^*W\,=\, E$.

Using the given condition $f^*W\,=\, E$, and the projection formula, we have
$$f_*E\,=\, f_*f^* W \,=\, W\otimes f_*{\mathcal O}_X.$$ Since ${\mathcal O}_Y\, \subset\,
f_*{\mathcal O}_X$ (see \eqref{e0}), we have
$$
W\, \subset\, W\otimes f_*{\mathcal O}_X\,=\, f_*E\, .
$$
Note that
$$
\mu(W)\,=\, \frac{1}{{\rm degree}(f)}\cdot \mu(E)\,.
$$

To prove the converse, assume that $f_*E$ contains a stable vector bundle $F$ such that
$$
\mu(F)\,=\, \frac{1}{{\rm degree}(f)}\cdot \mu(E)\,.
$$
Then $f^*F$ is stable by Theorem \ref{thm2}.
Let
$$
h\,:\, f^*F\, \longrightarrow\, E
$$
be the natural homomorphism; this $h$ is given by the inclusion map $F\, \hookrightarrow\, f_*E$
using the isomorphism in \eqref{eq_adjoint}. Since
\begin{enumerate}
\item both $f^*F$ and $E$ are stable,

\item $\mu(E)\,=\, {\rm degree}(f)\cdot \mu(F)\,=\, \mu(f^*F)$, and

\item $h\, \not=\, 0$,
\end{enumerate}
we conclude that $h$ is an isomorphism.
\end{proof}

\begin{remark}
Note that the proof of Theorem \ref{thm3} shows the following. Let $f\, \colon \, X\,
\longrightarrow\, Y$ be a finite surjective map of irreducible normal projective
varieties. Let $E$ be a stable vector bundle on $X$. If there is a vector bundle bundle
$W$ on $Y$ with $f^*W$ isomorphic to $E$, then $f_*E$ contains a stable vector bundle
$F$ such that $\mu(F)\,=\, \frac{1}{{\rm degree}(f)}\cdot \mu(E)$.
\end{remark}

\end{document}